\documentclass{amsart}
\usepackage{verbatim}
\usepackage{amssymb}
\usepackage{amscd}
\usepackage{amsmath}

\newtheorem{theorem}{Theorem}[section]
\newtheorem{case}{}
\newtheorem{corollary}{Corollary}[section]
\newtheorem{lemma}{Lemma}[section]
\newtheorem{proposition}{Proposition}[section]
\numberwithin{equation}{subsection}

\newtheorem{example}{Example}
\newcommand{\Q}{\mathbb Q}
\renewcommand{\H}{\mathbb H}
\newcommand{\G}{\mathcal G}
\newcommand{\Z}{\mathbb Z}
\newcommand{\p}{\mathfrak p}
\newcommand{\m}{\mathbf m}
\newcommand{\R}{\mathbb R}
\newcommand{\f}{\mathbf f}

\newcommand{\D}{\mathbf D}
\renewcommand{\c}{\mathbf c}
\renewcommand{\d}{\mathbf d}
\renewcommand{\m}{\mathbf m}
\renewcommand{\i}{\mathbf i}
\newcommand{\romanenum}{\renewcommand{\theenumi}{\roman{enumi}}}
\DeclareMathOperator{\sgn}{sgn}

\begin{document}

\bibliographystyle{plain}
\title{
Hilbert modular forms of weight 1/2 and theta functions}
\author{Sever Achimescu}
 \author{Abhishek Saha}
\date{March 2008}
\address{Institute of Mathematics "Simion Stoilow" of the Romanian Academy\\
P.O. Box 1-764\\
RO --  014700, Bucharest\\
Romania}
\address{Department of Mathematics 253-37 \\ California Institute of Technology \\ Pasadena, California 91125 \\ USA}

\email{achimesc@imar.ro}
\email{saha@caltech.edu}
\thanks{The authors wish to thank their advisor D.
Ramakrishnan for guidance and helpful discussions} \maketitle
\begin{abstract}
Serre and Stark found a basis for the space of modular forms of
weight $1/2$ in terms of theta series. In this paper, we generalize
their result - under certain mild restrictions on
the level and character - to the case of weight $1/2$ Hilbert modular forms over a totally
real field of narrow class number 1. The methods broadly follow those of Serre-Stark; however we are forced to overcome technical difficulties which arise when we move out of $\Q$.
\end{abstract}

\section*{Introduction}

\subsection{The problem and the main result} Shimura, at the end of his fundamental paper \cite{Shimura3} on elliptic
modular forms of half-integral weight,  mentioned certain
questions that were open at the time: one of them asked whether every modular form of weight 1/2 is a
linear combination of theta series in one variable. This was
answered in the affirmative by Serre--Stark \cite{Serre-Stark} who
gave an explicit basis for the space of modular forms of weight 1/2,
level $N$ and character $\psi$ in terms of certain theta series.
These theta series are denoted by $\theta_{\chi, t}$ where $\chi$ is
a primitive Dirichlet character and $t$ a positive integer so that
$\chi$ and $t$ are related in a precise manner to $N$ and $\psi$.
Such an explicit result has several nice applications, see for
instance Tunnell's work~\cite{Tunnell} on the ancient congruent
number problem.

It seems natural to generalize the Serre-Stark theorem to fields other than $\Q$, that is, to find an explicit basis in terms of theta
series for \emph{Hilbert modular forms} of weight $( \frac{1}{2},
\frac{1}{2}, ... \frac{1}{2})$. In this paper we achieve that in the case of a
totally real field $F$ of narrow class number 1 when the level
$\c$ and character $\psi$ of the form have certain nice properties.
In particular we assume that no prime dividing $\c$ splits in the
extension $F/\Q$ and that the Dirichlet character $\psi$ of the form
is trivial at the units (or equivalently, the corresponding finite order Hecke
character is trivial at all infinite places).
Under these assumptions we prove that the space of Hilbert modular forms of
weight $( \frac{1}{2}, \frac{1}{2}, ... \frac{1}{2})$, level $\c$
and character $\psi$ has a basis consisting of theta series that are almost
identical to the ones in Serre-Stark's theorem.

We note here that Shimura proved
(see Theorem~\ref{t:shimuragen}) that the space of weight $1/2$ Hilbert modular forms of \emph{all levels}
is \emph{spanned} by certain theta series; however his results do not seem
to give a \emph{basis}, nor do they appear to apply to a \emph{particular level}. Also, as noted by Deligne in a letter (appended at the end of \cite{Serre-Stark}) the problem can be attacked using the tools of representation theory. This was carried through successfully by Gelbart--Piatetski-Shapiro~\cite{gelpia}; however their result, like Shimura's, only finds a spanning set and also does not consider the levels.

We now briefly state the main result. Let $F$ be a totally real number field of narrow class number 1 and degree $n$ over $\Q$. Let $R$ be its ring of integers, $R^+ \subset R$ the subset of totally positive elements and $U$ the subgroup of units in $R$. For an ideal $\c$ of $R$ that is divisible by 4 and all of whose prime divisors are non-split, and a primitive Dirichlet character $\psi$ trivial on $U$, we let $M(\c, \psi)$ denote the space of Hilbert modular forms over $F$ of parallel weight $1/2$, level $\c$ and character $\psi$. For a primitive Dirichlet character $\chi$ trivial on $U$ and of conductor $r(\chi)$, and an element $t \in R^+$, we define the theta-series $\theta_{\chi, t}$ on the $n$-fold product of the upper half plane by $$\theta_{\chi, t}(z)= \sum_{x\in R} \chi^{-1}(x)e^{\pi i tr(x^2z)}.$$ Then our main theorem says the following.

\emph{A basis for $M(\c, \psi)$ is obtained by taking all the theta-series $\theta_{\chi, t}$ where we let $t$ vary over a set of representatives of $R^+/U^2$ and let $\chi$ satisfy, in addition to the conditions mentioned above, the following:
\begin{enumerate}
\item $4 r(\chi)^2t$ divides $\c$,
\item $\psi= \chi \epsilon_t$ where $\epsilon_t$ is the character associated to the quadratic extension $F(\sqrt{t})$.
\end{enumerate}    }

A few words about our methods. Though our techniques are similar to those of~\cite{Serre-Stark},
there are certain complications which arise because we are no longer
dealing with $\Q$; as a result many of the proofs of~\cite{Serre-Stark} do not extend to our case easily. Here are two main points of difference:

First, in section~\ref{s:oper} we prove various properties of certain operators (such as the symmetry operator) that are crucial to the theory of newforms for Hilbert modular forms of parallel weight $1/2$. In \cite
{Serre-Stark} these properties can be easily checked by hand and are left as exercises; however that is not the case here because we \emph{do not} have a simple closed formula for the automorphy factor. So we use an expression for the automorphy factor from Garrett's book~\cite{Garrett} and certain relations due to Shimura (and do some messy computations) to prove these properties. Furthermore we have to be very careful in the way we normalize these operators (and take into consideration the fact that the different of the field is no longer equal to 1) so that things work out.

Secondly, the proof of the crucial Theorem~\ref{t:newtheta} does not quite go through in a manner similar to~\cite{Serre-Stark}; the clever divisibility argument at the end of that proof breaks down here because of primes above $2$. We use a completely different method to get around this conundrum; we essentially use the fact that the size of the Fourier coefficients is bounded by Shimura's work.

Thus, the basis problem for the Hilbert modular case is not a
completely straightforward extension of \cite{Serre-Stark}, which,
we hope, justifies this article. Besides, we indicate, in a short section at the end, a motivation for solving this problem, by pointing out two potential applications which we hope to take up elsewhere.

\subsection{Structure of the paper}
In Section~\ref{s:prelim} we lay down notation, give some important
definitions and results that will be used throughout the paper,
state an important result due to Shimura and give the precise
statement of our main theorem.

In Section~\ref{s:hecke} we define the Hecke operators and write down
their action on Fourier coefficients.

Section~\ref{s:easy}, titled `Easy pickings' is the analogue of
\cite[Section 5]{Serre-Stark}. All the proofs carry over
\textit{mutatis mutandis} from there. We have included them for
completeness.

Section~\ref{s:oper} is similarly analogous to \cite[Subsection 3.4]{Serre-Stark}. In this section we define some important operators
(there are some differences from the corresponding definitions in
\cite{Serre-Stark} which arise because our definition of a modular
form is not \emph{quite} the same as Serre-Stark's) and prove the
same results as in there. However the calculations now are of a
higher order of difficulty than in \cite{Serre-Stark} because, unlike
in the classical case, there is no simple formula for the automorphy
factor. As a result the proofs are more technical. This is probably
the hardest part of the paper, involving messy computations.

In Section~\ref{s:newforms} we outline the theory of newforms for
our purposes. The proofs are but formal consequences of the results
of the previous two sections and essentially identical to the
corresponding proofs in \cite{Serre-Stark}. Therefore we do not
include them.

In Section~\ref{main}, we define the $L$-series and use it to
characterize a newform. Using that, we prove our main
theorem. At the end of this section, we illustrate our theorem by writing down bases for the spaces of weight $1/2$ Hilbert modular forms over  $\Q(\sqrt{2})$ for various levels.

Finally, in Section~\ref{s:app} we mention some potential applications of our work.

\section{Preliminaries}\label{s:prelim}
\subsection{Notation} Let
$F$ be a totally real number field, $R$ its ring of integers, $D$ its discriminant and
$\delta$ its different. By abuse of notation we also use $\delta$ to
denote a fixed totally positive generator of the different. We
assume that $F$ has narrow class number one and we let $n$ denote
the degree of $F/\Q$. Let the group of units of $F$ be denoted by
$U$ and the group of totally positive units by $U^2$ (since the
field is of narrow class number one, all totally positive units are
squares). For any $t \in F$, we use the notation $t
>> 0$ to mean that $t$ is totally positive.

We denote the adelization of $F$ by $F_\mathbf{A}$ and the ideles by
$F_\mathbf{A}^{\times}$. For any $x \in F$ let $N(x)$ denote its
norm over $\Q$. For an ideal $\m \subset R$, we will let $N(\m)$
denote the cardinality of $R/\m$. Let $\infty$ denote the set of
Archimedean places of $F$ and $\mathbf{f}$ denote the finite places. For $g \in F_\mathbf{A}^{\times}$ we denote $g_{\m} = \prod_{v | \m}g_v$ and $g_{\infty}=\prod_{v \in \infty}g_v$.
For $v \in \infty$ we denote the positive elements of $F_v$ by
$F_v^{+}$. By $F_{\infty}^\circ$ we mean the connected component at
infinity of the identity, i.e. $$F_{\infty}^\circ= \prod_{v \in
\infty } F_v^{+} \simeq (\R^+)^n.$$

Let $\mathbb{H}^n$ (resp. $\mathbb{C}^n$) denote the n-fold product
of the upper half plane (resp. complex plane). For $z =
(z_1,..,z_n)$ in $ \mathbb{C}^n$ or $\R^n$ and any $\alpha \in \R$ we
put $$ z^\alpha = \prod_{i=1}^n z_i^\alpha,\quad e(z) =
\prod_{i=1}^n e^{2 \pi i z_i},\quad N(z) = \prod_{i=1}^n z_i, \quad
tr(z)= \sum _{i=1}^n z_i.$$ We also use the symbol $e(z)$ for $z \in
F$ using the $n$ embeddings of $F$ in $\R$.

Furthermore, for any prime(i.e. a finite place) $p$ in $\R$ we
define the character $e_p$ on $F$ as follows: For $x \in F$ let
$$e_p(x) = e^{-2\pi i y}$$ where $y\in \cap_{q \neq p'}(\Z_q \cap
\Q)$, $y- Tr_{F_p/\Q_{p'}}(x) \in \Z_{p'}$. Here $q$ is any prime in
$\Z$ and $p'$ is the prime below $p$, that is, $p' = p \cap \Q.$

\subsection{Conventions on characters}\label{s:dirichlet}
Let $\m$ be an ideal of $R$. A Dirichlet character mod $\m$ is a
function $\phi$ from $R$ to the unit circle such that:
\begin{enumerate}
\item There exists a homomorphism $\overline{\phi}$ from the finite
group $(R/\m)^\times$ to the unit circle such that for any $a$ in
$R$ that is relatively prime to $\m$ we have $\phi(a) =
\overline{\phi}(\overline{a})$

\item For any $a$ that shares a common factor with $\m$, we have
$\phi(a) = 0$.

\end{enumerate}
Such a $\phi$ is called \emph{primitive} if
$\overline{\phi}$ does not factor through $(R/\m')^\times$ for some
proper divisor $\m'$ of $\m$.

By a Hecke character of $F$ we mean a character of
$F_\mathbf{A}^{\times}$ which is trivial on $F^{\times}$ and has
values in the unit circle. For a Hecke character $\psi$ and any
place $v$, $\psi_v$ denotes the restriction of $\psi$ to $F_v^\times$. For
an ideal $\c$ of $R$, let $$\psi_\c = \prod_{v \mid \c} \psi_v,$$  $$\psi_\f = \prod_{v | f}\psi_v$$
and
$$\psi_\infty = \prod_{v \in \infty} \psi_v.$$ The conductor of $\psi$
always refers to its finite part and is denoted by $r(\psi)$. In a mild abuse of notation we will use, for $g$ in $F^\times$ or even $F_\mathbf{A}^{\times}$, $\psi_\c(g)$ (resp. $\psi_\infty(g)$) to mean $\psi_\c(g_\c)$ (resp. $\psi_\infty(g_\infty)$). Also
let $\psi^\ast$ denote the corresponding character on the ideals as
defined in \cite[p. 238]{Shimura4}. In particular, if $I$ is an
ideal of $R$ generated by $s$ and $(I , r(\psi)) = 1$, we have
$$\psi^\ast(I) = \overline{\psi_\c(s)\psi_\infty(s)}$$ for any $\c$
divisible by $r(\psi)$ with $(I, \c)=1$. We will use the notation $\psi^\ast(a)$ for
$a\in R$ to denote $\psi^\ast((a))$.

For any $\tau \in F$ let
$\epsilon_\tau$ denote the Hecke character of $F$ corresponding to
$F(\tau^{1/2})/F$.

\textbf{Comment}: It is well known that any \emph{finite order}
Hecke character $\psi$ of $F$ gives rise to a primitive Dirichlet character
mod $r(\psi)$. This correspondence is bijective. Moreover, for such a finite order Hecke character $\psi$, we have $\psi_\infty(g) = \prod_{v \in \infty}\sgn(g_v)^{e_v}$ with each $e_v = 0$ or $1$. Thus $\psi_\infty$  is trivial if and only if each $e_v=0.$ It can be checked that this happens if and only if the corresponding Dirichlet
character is trivial on the units of $R$.

\subsection{Conventions on modular forms}

\emph{In the rest of this paper, unless mentioned otherwise, we will use the term Hecke character to mean finite order Hecke character.}

Given a $2 \times 2$ matrix $\alpha= \left( \begin{array}{ccc}
a & b  \\
c & d  \\
\end{array} \right)$ we write $a=a_\alpha, \ b=b_\alpha, \
c=c_\alpha$ and $d= d_\alpha$.

Let $G = SL_{2}(F)$. We will consider weight
$(\frac{1}{2},...,\frac{1}{2})$ modular forms on the congruence
subgroups of $G$.

For any two fractional ideals $\mathfrak{f}, \mathfrak{g}$ of $R$,
let $\Gamma[\mathfrak{f}, \mathfrak{g}]$ denote the subgroup of $G$
consisting of matrices $\gamma$ such that
 $a_{\gamma}, \ d_{\gamma} \in R$, $b_{\gamma} \in \mathfrak{f}, \
c_{\gamma} \in \mathfrak{g}$. Let $\D$ denote the group $\Gamma[2
\delta^{-1}, 2 \delta]$. A congruence subgroup is a subgroup of $\D$
that contains a principal congruence subgroup $\Gamma(N)$ for some
integer $N$, where
       $$\Gamma(N) = \{\gamma \in SL_2(R) :  \gamma \equiv I \pmod{N} \}.$$

For any two integral ideals $\c,\d$, with $4 | \c$ we use the
notation $$\Gamma_{\c,\d} = \Gamma[2\delta^{-1}\d,2^{-1}\delta\c]$$
and $$\Gamma_\c = \Gamma[2\delta^{-1},2^{-1}\delta\c]. $$ Note that
$\Gamma_{\c,\d}$ and $\Gamma_\c$ are congruence subgroups.

For $\gamma \in \D$, $z \in \mathbb{H}^n$, let $h(\gamma, z)$ denote
the automorphy factor $\frac{\theta(\gamma z)}{\theta(z)}$, where
$$\theta(z) =\sum_{x \in R} e(x^2 z /2).$$  A generalization of this
automorphy factor is introduced in \cite{Shimura2} where many of its
properties are proved.

Now, let $\gamma \in \D$, and $f$ be a holomorphic function on
$\H^n$. We use the notation $$(f
\parallel \gamma)(z) = h(\gamma, z)^{-1}f(\gamma(z)).$$

Suppose now that $\c$ is an ideal as above and $\psi$ is a  Hecke
character whose conductor divides $\mathbf{c}$ and $\psi_\infty(-1)
= 1$. Let $M(\mathbf{c} , \psi)$ denote the space of modular
forms of weight $( \frac{1}{2},  ... \frac{1}{2})$ on
$\Gamma_\mathbf{c}$ with character $\psi$. In other words,
$M(\mathbf{c} , \psi)$ is the set of holomorphic functions $f$ on
$\H^n$ satisfying $$f
\parallel \gamma = \psi_\mathbf{c}(a_\gamma)f$$
for all $\gamma \in
\Gamma_\mathbf{c}$. Note that our definition follows \cite{Shimura1}
(and is slightly different from \cite{Serre-Stark} where $f$
satisfies $f
\parallel \gamma = \psi_\mathbf{c}(d_\gamma)f$).

For each such $\mathbf{c}$ let $M^1(\mathbf{c})$ be the union of all
$M(\mathbf{c} , \psi)$ with $\psi$ varying over all  Hecke characters
with conductor dividing $\mathbf{c}$ and $\psi_\infty(-1) = 1.$
Let $M^1$ be the union of all the $M^1(\mathbf{c})$ as $\mathbf{c}$
varies over the integral ideals of $R$ divisible by 4. Finally, let
$M$ be the space of all weight $( \frac{1}{2},  ... \frac{1}{2})$ modular forms on congruence subgroups of
$G$. Clearly for any such $\mathbf{c} , \psi$, $M(\mathbf{c} , \psi)
\subset M^1(\mathbf{c}) \subset M^1 \subset M$.

Any $f \in M$ has a \emph{Fourier expansion}
$$f(z) = \sum_{\xi \in F} a(\xi)e(\xi z /2).$$ We call $a(\xi)$ the
Fourier coefficient for the place $\xi$.

If $f$ belongs to $M^1$ then by \cite[p.~780]{Shimura1}, the
Fourier coefficients associated to places outside $R$ are zero. Thus
$f$ has a Fourier expansion
$$f(z) = \sum_{\xi \in R} a(\xi)e(\xi z /2).$$

We are interested in the question of finding a basis for each of the
spaces $M(\mathbf{c} , \psi)$.

\subsection{Theta functions}

Let $\eta$ be a locally constant function on $F$, i.e. a
complex valued function for which there exists two
$\mathbb{Z}$-lattices $L$ and $M$ in $F$ such that $\eta(x)=0$ for
$x$ not in $L$ and $\eta(x)$ depends only on $x$ modulo $M$.

The following alternate criterion will be useful.

\begin{proposition}\label{p:criterion}
A function $\eta : F \rightarrow \mathbb{C}$ is locally constant if
and only if there exist integers $m,n$ such that $\eta(x) = 0$ for
$x$ not in $\frac{1}{m} R$ and $\eta(x)$ depends only on $x$ mod
$(n)$.
\end{proposition}

\begin{proof} Any $\mathbb{Z}$-lattice contains $(n)$ and is contained in
$(\frac{1}{m})$ for some $m,n$.
\end{proof}

Let $\mathfrak{L}(F)$ denote this space of locally constant
functions. We define the function $\theta_{\eta}$ on $\H ^n$ by
$$\theta_{\eta}(z) = \sum_{\xi \in F} \eta( \xi) e(\xi^2 z
/2).$$

We have the following proposition.

\begin{proposition}
Let $\eta \in \mathfrak{L}(F)$. Then $\theta_{\eta} \in M$.
\end{proposition}

\begin{proof}
 This follows from \cite[Lemma 4.1]{Shimura1}. Indeed the
proof there makes it clear that $\theta_\eta$ is a modular form for
the largest congruence group contained in $\{\alpha \in D, ^\alpha
\eta = \eta\}$, where $^\alpha \eta$ denotes the action of $\alpha$
on $\eta$ as described in \cite[p. 775]{Shimura1}.

\end{proof}

\subsection{An important example}\label{s:exam}
The following example from \cite{Shimura1} introduces the theta
series that is fundamental to this paper.

\begin{example} [\cite{Shimura1} , p. 784-785]

Let $\chi$ be a  Hecke character of $F$ of conductor $f$ such that
$\chi_\infty(-1) = 1$. Suppose $\omega_v$ denotes the characteristic
function of $R_v$ and let $$\eta(x) = \prod_{v\in \mathbf{f}}
\eta_v(x_v)$$ where:

\begin{itemize}

\item $\eta_v = \omega_v$ if $v \nmid f$

\item $\eta_v = \chi_v(t)^{-1}$ if $v \mid f$ and $\mid t\mid _v = 1$

\item $\eta_v = 0$ if $v \mid f$ and $\mid t\mid _v \neq 1$

\end{itemize}
Then $\theta_\eta(z) \in M(4f^2 , \chi)$.
\end{example}

For any Hecke character $\chi$ of $F$  such that $\chi_\infty(-1) =
1$ we define $\theta_\chi$ to equal $\theta_\eta$ where $\eta$ is as
in the above example. Thus $\theta_\chi \in M(4r(\chi)^2 , \chi).$

For any totally positive $t \in F$ let $\theta_{\chi,t}(z) :=
\theta_\chi(tz)$. We have $\theta_{\chi, t} \in M(\c , \psi)$
whenever $(4r(\psi)^2t) \mid \c$ and $\psi = \chi\epsilon_t$. Refer
to Lemma~\ref{l:shift} for a proof of this fact. Similarly, for any
function $\eta \in \mathfrak{L}(F)$ let $\theta_{\eta,t}(z):=
\theta_\eta(tz)$.

\subsection{Two generating sets}

The following important theorem is due to Shimura and is contained
in \cite{Shimura1}.

\begin{theorem}[Shimura] \label{t:shimuragen}
$M$ is spanned by the functions $\theta_{\eta,t}$ for $t \in F$
totally positive and $\eta \in \mathfrak{L}(F).$
\end{theorem}

What about the space of forms $M^1$?

We make the following preliminary observations:

Any $f \in M$, by the above theorem, can be written as
\begin{equation}\label{e:eqftheta}
f(z) = \theta_{\eta_1}(t_1z) + \theta_{\eta_2}(t_2z) + ..+
\theta_{\eta_k}(t_k z)
\end{equation}

with $0 <<t_i \in F$.

Replacing each $\eta_i(z)$ by $\frac{\eta_i(z) + \eta_i(-z)}{2}$, we
may assume that $\eta_i(z) = \eta_i(-z)$. Note that this does not
change the functions $\theta_{\eta_i}$.

Also, we may assume that the $t_i$ are distinct mod $(F^{\ast})^2$.
For, if $t_1 = s^2t_2$, say, then $$\theta_{\eta_1}(t_1z) +
\theta_{\eta_2}(t_2z) = \theta_\eta(t_2z)$$ where $\eta(z) =
\eta_1(z/s) + \eta_2(z)$, and so we may combine those two summands
into a single one.

Furthermore, if $\eta_i(\xi) = 0$ for $\xi$ not in $(\frac{1}{m})R$,
then $\eta_i(\frac{\xi}{m}) = 0$ for $\xi$ not in $R$. Moreover,
observe that $\theta_{\eta_i}(t_iz) = \theta_{\eta_i^{'}}(t_iz/m^2)$
where $\eta_i^{'}(z) = \eta_i(z/m)$. So, in~\eqref{e:eqftheta} we may assume that
each $\eta_i$ is 0 outside $R$. We can now give a set of generators
for $M^1$.
\begin{theorem}\label{t:shimuragen2} $M^1$ is spanned by the functions $\theta_{\eta,t}$ where $t \in R$ is totally positive,
and  $\eta \in \mathfrak{L}(F)$ satisfies $\eta(z) =0$ if $z$ does
not belong to $R$.
\end{theorem}

\begin{proof} Any $f \in M^1$, by the above comments can be written as

\begin{equation}\label{e:ftheta2} f(z) = \theta_{\eta_1}(t_1z) + \theta_{\eta_2}(t_2z) + ..+
\theta_{\eta_k}(t_k z)
\end{equation}
where $0 <<t_i \in F$ are distinct mod $(F^{\ast})^2$ and $\eta_i
\in \mathfrak{L}(F)$ are $0$ outside $R$.

Then, because the $t_i$ are distinct mod $(F^{\ast})^2$ the various
$\theta_{\eta_i}(t_iz)$ contribute distinct terms to the Fourier
expansion of $f$. However only the Fourier coefficients
corresponding to elements of $R$ can be non-zero.

So for each $i$ we must have $\eta_i(\xi) = 0$ whenever $\xi^2t_i$
not in $R$. For a fixed $t_i$, the set of $\xi \in R$ such that
$\xi^2t_i \in R$ is an ideal, hence generated by some $h$. Put
$\eta_i^{'}(z) = \eta_i(hz)$. Then $\theta_{\eta_i}(t_iz) =
\theta_{\eta_i^{'}}(t_i h^2z)$. Thus replacing $\eta_i$ by
$\eta_i^{'}$ and $t_i$ by $t_i h^2$ we see that $\eta_i^{'}$ is
still 0 outside $R$, but now $t_i h^2$ also belongs to $R$.

In other words we have shown that in~\ref{e:ftheta2}, under the assumption $f
\in M^1$ we can have $0<<t \in R$, and $\eta$ is $0$ outside $R$.

Conversely any such sum is in $M^1$ by \cite[Proposition
3.2]{Shimura1} and \cite[p. 154]{Garrett}.

This completes the proof. \end{proof}

\begin{corollary} \label{c:bound}
Let $f(z) = \sum_{\xi \in R} a(\xi)e(\xi z /2)$ be an element of
$M(\c, \psi)$ for some $(\c, \psi)$. Then there is a constant $C_f$
such that $\mid a(\xi)\mid < C_f$ for all $\xi \in R$.
\end{corollary}
\begin{proof}
By the above theorem, it suffices to prove that $\theta_{\eta,t}$
has this property. But that follows easily from
Proposition~\ref{p:criterion}.
\end{proof}

\subsection{Statement of the main theorem}\label{s:main}
Let $R^+$ denotes the set of
totally positive elements in $R$. Fix a complete set of
representatives $T$ of $R^{+}/U^2$.

Suppose $\c$ is an integral ideal and $\psi$ a
 Hecke character. Define $\Omega(\c,
\psi)$ to be the set of pairs $(\chi, t)$ such that:
\begin{enumerate}
\item $\chi$ is a  Hecke character with $\chi_\infty$ trivial and $t \in T.$
\item $4r(\chi)^2t$ divides $\c.$
\item $\psi = \chi\epsilon_t.$
\end{enumerate}

Recall the definition of $\theta_{\chi, t}$ from
Section~\ref{s:exam}. Our main Theorem is as follows:
\begin{theorem} \label{t:main}Suppose $\c$ is an integral ideal divisible by $4$. Let $\psi$ be a
 Hecke character of $F$ such that $\psi_\infty$ is
trivial and $r(\psi)$ divides $\c$. Assume that any prime ideal $\p$
dividing $\c$ has the property that $\p$ is the unique prime
ideal of $R$ that lies above $\p \cap \Z$. Then the functions $\theta_{\chi,t}$ with $(\chi, t) \in
\Omega(\c, \psi)$ form a basis of $M(\c, \psi).$
\end{theorem}

We prove this theorem in section~\ref{main}.
\section{Hecke operators}\label{s:hecke}
\subsection{Some definitions}

Let $GL_2^+(F)$ denote the subgroup of $GL_2(F)$ consisting of
matrices whose determinant is totally positive. Let $\G$ denote the
group extension of $GL_2^+(F)$ consisting of pairs $[A, \phi(z) ]$
where $A =\left( \begin{array}{ccc}
a & b  \\
c & d  \\
\end{array} \right) \in GL_2^+(F)$ and $\phi(z)$ is a holomorphic
function on $\H^n$ satisfying $\phi(z)^2 = t N (\mbox{det} A)^{-1/2}
\prod (c^{(i)}z_{i} + d^{(i)}) $ where $A^{(i)} =\left(
\begin{array}{ccc}
a^{(i)} & b^{(i)}  \\
c^{(i)} & d^{(i)}  \\
\end{array} \right) $ are the various embeddings of $A$ in $GL_2(\R)$ and $t$ is a complex number with $\mid t \mid =
1$. The group law in $\G$ is given by $[A, \phi(z)][B, \psi(z)] =
[AB, \phi(Bz)\psi(z)]$.

The group $\G$ acts on the $\emph{right}$ of the space of
holomorphic functions on $\H^n$ as follows: For a holomorphic
function $f$ on $\H^n$ define $f \mid [A, \phi(z)] =
\phi(z)^{-1}f(Az)$. Note also that the group $D$ embeds in $\G$ via
$A \rightarrow [A, h(A, z)$. Furthermore, we have $(f \parallel
A)(z) = f \mid [A, h(A , z)]$.

For $\gamma = w_1tw_2$ where $w_1, w_2 \in \D$ and $t
=\begin{pmatrix}
1/a & 0  \\
0 & a  \\
\end{pmatrix}$ for some $a \in R$ , define $J_\Xi(\gamma, z)= h(w_1w_2,z)$ .  The quantities $J_\Xi(\gamma, z)$ and
$h(\gamma, z)$ coincide whenever $\gamma \in \D$. We also recall from~\cite{Shimura1} that:

\begin{enumerate}

\item $J_\Xi\left(\begin{pmatrix}
1/p & 2b/(\delta p)  \\
0 & p  \\
\end{pmatrix}, z\right) = N(p)^{1/2}$ for a prime $p$ and element $b$ in $R$.

\item $J_\Xi\left(\begin{pmatrix}
1 & 2h/(\delta p)  \\
0 & 1 \\
\end{pmatrix}, z\right) = N(p)^{1/2} \left(\sum_{x\in
(R/p)}e_p(hx^2/(p\delta) )\right)^{-1}$ where $p$ is a prime and $h
\in R$ is not divisible by $p$.

\item $J_\Xi\left(\begin{pmatrix}
p & 0  \\
0 & 1/p  \\
\end{pmatrix}, z\right) = N(p)^{-1/2}$.\\
\end{enumerate}

A key property of $J_\Xi$ is that it is a partial automorphy factor.
To be precise, it has the following properties (see \cite{Shimura1}):\\

(a)$J_\Xi(y_1xy_2 , z) = h(y_1, xy_2(z))J_\Xi(x, y_2(z))h(y_2, z)$,
if $y_1 , y_2$ belong to $\D$.

(b)$J_\Xi(k^{-1}, z) = J_\Xi(k, k^{-1}(z))^{-1}$, where $k \in \D \sigma
\D$ with $\sigma =\begin{pmatrix}
1/a & 0  \\
0 & a  \\
\end{pmatrix}$ with $a$ relatively prime to $2$.

Let $\gamma \in \Gamma_{\c}.$ We now give a complicated, but still useful, formula for $h(\gamma,z)$.

For $d \in R - \{0\}$, define $$\epsilon(d) = (i \sgn d)^{1/2}2^{-n/2} D^{-1/2} \sum_{v\in \delta^{-1}/2R}e(-v^2d/4).$$ We also define $\widetilde{\epsilon}(d) = i^s$ where $s$ is the number of negative embeddings of $d$. Then~\cite[p. 142]{Garrett} tells us that \begin{equation}\label{e:hformula}h(\gamma, z) =  \epsilon(d_\gamma) \widetilde{\epsilon}(d_\gamma) (\epsilon_{c_\gamma})^\ast(a_\gamma)(c_\gamma z + d_\gamma)^{1/2}\end{equation}

Also, by~\cite[p. 146]{Garrett} we have $$\theta\left(\begin{pmatrix}0& -\delta^{-1}\\ \delta &0\end{pmatrix} z\right) / \theta(z) = (- \i z)^{1/2}N(\delta)^{1/2}.$$  We extend the notation $h(\gamma, z)$ to this case by defining \begin{equation}\label{hwformula}h\left(\begin{pmatrix}0& -\delta^{-1}\\ \delta &0\end{pmatrix}, z\right) = (- \i z)^{1/2}N(\delta)^{1/2}\end{equation}

For each totally positive prime element $p \in R$ we define a
\emph{Hecke operator} $T_{p^2}$ on $M(\c, \psi)$ that sends $f $ to
$f \mid T_{p^2}$ where

\begin{align*}
f \mid T_{p^2} = N(p)^{-3/2}&\overline{\psi_\infty(p)}\quad \bigg(\sum_{b
\in R/p^2} f \mid \left[
\begin{pmatrix}
1/p & 2b/(\delta p)  \\
0 & p  \\
\end{pmatrix},J_\Xi\left(\begin{pmatrix}
1/p & 2b/(\delta p)  \\
0 & p  \\
\end{pmatrix}, z\right)   \right] \\
&+ \quad \overline{\psi_\c(p)}\sum_{h\in (R/p)^\times} f  \mid
\left[
\begin{pmatrix}
1 & 2h/(\delta p)  \\
0 & 1  \\
\end{pmatrix},   J_\Xi
\left(\begin{pmatrix}
1 & 2h/\delta p  \\
0 & 1  \\
\end{pmatrix} ,z \right)  \right] \\
&+\quad  \overline{\psi_\c(p^2)} f  \mid \left[
\begin{pmatrix}
p & 0  \\
0 & 1/p \\
\end{pmatrix}, J_\Xi\left(\begin{pmatrix}
p & 0  \\
0 & 1/p  \\
\end{pmatrix}, z\right)   \right]\quad \bigg).
\end{align*}

\subsection{Action of the Hecke operator on Fourier coefficients}
The next proposition, which is a restatement of (\cite{Shimura1},
Proposition 5.4) gives the explicit action of $T_{p^2}$ on the
Fourier coefficients of a modular form. In particular it also shows
that if $p$ and $p'$ are two totally positive elements that generate
the same prime ideal, then $T_{p^2}$ coincides with $T_{p'^2}$

\begin{proposition}[Shimura]\label{p:Hecke}
Suppose $p$ be a totally positive prime element of $R$ and $f \in
M(\mathbf{c}, \psi)$ is given by
$$f(z) = \sum_{\xi \in R} a(\xi)e(\xi z /2).$$

Then
$$(f \mid T_{p^2})(z) = \sum_{\xi \in R}
b(\xi)e(\xi z /2)$$

where $$\psi_\infty(p)b(\xi) = a(\xi p^2) +  \left\{
\begin{array}{ll}
\overline{\psi_\mathbf{c}(p)}N(p)^{-1}(\frac{\xi}{p})a(\xi)
+ \overline{\psi_\mathbf{c}(p^2)}N(p)^{-1}a(\xi/p^2)& \textrm{if $p \nmid \mathbf{c}$}\\
\\
0 & \textrm{if $p \mid \mathbf{c}$}
\end{array} \right.$$

where $a(\xi/p^2) := 0$ if $p^2 \nmid \xi$.

\end{proposition}

\begin{corollary}\label{cbound2}

Suppose $$f(z) = \sum_{\xi \in R} a(\xi)e(\xi z /2)$$ is an element
of $M(\c, \psi)$ and $f\mid T_{p^2} = c_p f$ for some prime $p \mid
\c$. Then
$$a(\xi p^{2n}) = (\psi_\infty (p))^n c_p^n a(\xi)$$ and $\mid c_p
\mid \leq 1$.
\end{corollary}

\begin{proof} The assertion about $a(\xi p^2)$ follows from the
above Proposition. Now Corollary~\ref{c:bound} implies that $\mid
c_p \mid \leq 1$.

\end{proof}

\section{Easy pickings}\label{s:easy}
\subsection{Eigenforms of Hecke operators}
Consider the Petersson scalar product $<f,g>$ on $M(\mathbf{c},
\psi)$. The definition is analogous to the classical case, see for
instance \cite{Shimura1}.

By standard calculations~\cite[Proposition 5.3]{Shimura1},
$\overline{\psi^\ast(p^2)}T_{p^2}$ is a Hermitian operator if $p$
does not divide $\c$. Hence:

\begin{lemma}\label{l:basis} There is a basis of $M(\mathbf{c}, \psi)$ consisting of
eigenforms for all the $T_{p^2}$ where $p>>0$ is a prime in $R$ and
$p \nmid \mathbf{c}$.
\end{lemma}

So it is important to study the modular forms that are eigenvalues
for the Hecke operators. But first we prove an auxillary lemma.

\begin{lemma}\label{l:numbasis} The following hold:

(a) There is a basis of $M(\mathbf{c}, \psi)$ consisting
of forms whose coefficients belong to a number field.

(b)If $f(z) = \sum_{\xi \in R}a(\xi)e(\xi z /2) \in M(\mathbf{c},
\psi)$ has each $a(\xi)$ algebraic, then the $a(\xi)$ have bounded
denominators (i.e. there exists a non zero integer $D$ such that
$Da(\xi)$ is an algebraic integer for all $\xi$).\end{lemma}

\begin{proof}
(a) is just~\cite[Proposition 8.5]{Shimura1} while (b) follows
easily from Theorem~\ref{t:shimuragen2} above. \end{proof}

\begin{lemma} Let $f(z) = \sum_{\xi \in R}a(\xi)e(\xi z
/2) \in M(\mathbf{c}, \psi)$ be an eigenvector of $T_{p^2}$
with eigenvalue $c_p$ where $p \nmid \mathbf{c}$.
Suppose $0<< m \in R$ such that $p^2 \nmid m$. Then:

(a) $a(mp^{2n}) =
a(m)\overline{\psi_\mathbf{c}(p)^n}(\frac{m}{p})^n$ for every $n
\geq 0$

(b) If $a(m)\neq 0$, then $p \nmid m$ and $c_p
=\psi^\ast(p)(\frac{m}{p})(1 + N(p)^{-1})$
\end{lemma}

\begin{proof} Since $T_{p^2}$ maps forms with algebraic coefficients
into themselves, it follows from Lemma~\ref{l:numbasis} by simple
linear algebra that the eigenvalue $c_p$ is algebraic and that the
corresponding eigenspace is generated by forms with algebraic
coefficients. So we assume that the coefficients $a(\xi)$ are
algebraic.

Consider the power series $A(T)= \sum_{n=0}^{\infty}a(mp^{2n})T^n$

Using Proposition \ref{p:Hecke}, we get, by the same argument as in~\cite[p. 452]{Shimura3}.

$$ A(T) = a(m)\frac{1 - \alpha T}{(1 - \beta T)(1 - \gamma
T)}$$

where

$$\alpha =\overline{\psi_\mathbf{c}(p)}N(p)^{-1}(\frac{m}{p})$$

and

$$\beta + \gamma = \psi_\infty (p)c_p,\quad \beta\gamma =\overline{\psi_\mathbf{c}(p^2)}N(p)^{-1}.$$

This already implies that $a(m) = 0$ implies $a(mp^{2n}) =0 \forall
n$. Hence we may assume that $a(m) \neq 0$ in which case $A(T)$ is a
non zero rational function of $T$. Viewing $A(T)$ as a function over
a suitable finite extension of $\mathbb{Q}_p$, we see using
Lemma~\ref{l:numbasis}(b) that $A(T)$ converges in the $p$-adic unit
disk $U$ defined by $\mid T \mid_p < 1$; hence $A(T)$ cannot have a
pole in $U$. However, since $\beta \gamma =
\overline{\psi_\mathbf{c}(p^2)}N(p)^{-1}$ one of $\beta^{-1} ,
\gamma^{-1}$ belongs to $U$. Assume it is $\beta^{-1}$. Since $A(T)$
is holomorphic we must then have $\alpha = \beta$. So
$A(T)=\frac{a(m)}{(1 - \gamma T)}$ and so $a(mp^{2n}) = \gamma^n
a(m)$. Since $\beta \gamma \neq 0$ we have $\alpha \neq 0$, hence $p
\nmid m$. Moreover $\gamma = \beta \gamma / \alpha =
\overline{\psi_\mathbf{c}(p)}(\frac{m}{p})$. So $a(mp^{2n}) =
\gamma^n a(m) = a(m)\overline{\psi_\mathbf{c}(p)^n}(\frac{m}{p})^n$.
This proves (a) while (b) follows from $c_p =
\overline{\psi_\infty(p)}(\alpha + \gamma)$. \end{proof}

An element $t \in R$ is called squarefree if it is not divisible by
the square of a prime element of $R$.

\begin{theorem}\label{t:serre}
Let $$f(z) = \sum_{\xi \in R} a(\xi) e(\xi z / 2)$$ be a non zero
element of $M(\c, \psi)$ and let $\c'$ be an ideal of $R$ such that
$ \c \mid \c '$. Assume that for all primes $p \nmid \c'$ we have
$f\mid T_{p^2} = c_pf$ where $c_p \in \mathbb{C}$. Then there exists
a unique(up to multiplication by an unit) totally positive
squarefree element $t \in R$ such that $a(\xi) = 0$ unless
$\frac{\xi}{t}$ is the square of an element of $R$. Moreover

\begin{enumerate}

\item $t \mid \c'$

\item $c_p = \psi^\ast(p)(\frac{t}{p})(1 + N(p)^{-1})$ \ if $p \nmid
\c'$

\item $a(\xi u^2) = a(\xi)\overline{\psi_\mathbf{c}(u)}(\frac{t}{u})$
\ if $(u \ , \ \c') = 1$
\end{enumerate}

\end{theorem}

\begin{proof}. Let $\xi, \xi' \in R$ such that $a(\xi), a(\xi') \neq
0$. We first show that $\xi' / \xi$ is a square. Let $P$ be the set
of primes $p$ with $p \nmid (\c' \xi \xi')$. If $p\in P$, the
previous lemma shows that
$$\overline{\psi_\infty(p)}\overline{\psi_\mathbf{c}(p)}(\frac{\xi}{p})(1 + N(p)^{-1}) = c_p = \overline{\psi_\infty(p)}\overline{\psi_\mathbf{c}(p)}(\frac{\xi'}{p})(1 +
N(p)^{-1})$$

Hence $$(\frac{\xi}{p}) = (\frac{\xi'}{p})$$ for almost all p. But
this means that almost all primes split in the extension
$F(\sqrt{\xi\xi'}) / F$ and hence by a well known result the
extension must be trivial, i.e. $\xi' /\xi$ is a square. Write $\xi
= tv^2 , \xi' = tv'^2$ with $t$ totally positive square free.

This proves the first assertion of the theorem, i.e. the existence
of $t$. Now write $v= p^nu$ with $p \nmid \c'$ and $(p,u) = 1$. So $
\xi = t p^{2n} u^2$. Applying the previous lemma to $tu^2$ we have
$a(\xi) = a(tu^2)\overline{ \psi_ \c(p)^n }(\frac{tu^2}{p})^n$.
Hence $a(tu^2) \neq 0$ and part (b) of the lemma above shows that $p
\nmid t$ and $c_p = \psi^\ast(\p)(\frac{t}{p})(1 + N(p)^{-1})$.

Hence every prime factor of $t$ divides $\c'$ ; since  $t$ is
squarefree this implies $t \mid \c'$, and (1)~,(2) are proved. As
for (3), it is enough to check it for $u = p$ with $p \nmid \c'$ and
$u$ a unit. The case of $u = p$ follows from writing $\xi = \xi_0
p^{2a}$ with $p^2 \nmid \xi_0$ and applying part (a) of the previous
lemma, while the case of an unit follows from \cite{Shimura1},
Proposition 3.1. \end{proof}

\section{Some operators}\label{s:oper}

Note that all operators on spaces of modular forms defined in this
paper act on the \emph{right}. This is done so that composition of
operators is compatible with multiplication in $\G$.

Fix a totally positive generator $c$ of $\c$. We define the
following operators on $M(\c, \psi)$.

\begin{itemize}

\item(The shift operator) For any totally positive $m \in R$, the shift operator $V(m)$ is defined as
$$V(m) = N(m)^{-1/4}\left[\begin{pmatrix}
m & 0  \\
0  & 1 \\
\end{pmatrix}, N(m)^{-1/4}\right]$$

Thus $$(f \mid V(m))(z) = f(mz).$$

\item(The symmetry operator) The symmetry operator
$W(c)$ is defined as

$$W(c) = [W_0, J_\Xi(W_0,z)][V_0(c), N(c)^{-1/4}]$$ where $W_0
=\begin{pmatrix}
0 & -2\delta^{-1}  \\
2^{-1}\delta & 0  \\
\end{pmatrix}$ and $V_0(c) = \begin{pmatrix}
c & 0  \\
0 & 1  \\
\end{pmatrix}$.

Thus $$(f \mid W(c))(z) = f\mid \left[
\begin{pmatrix}
0 & -2\delta^{-1}  \\
2^{-1} \delta c & 0  \\
\end{pmatrix}, (- \i z)^{1/2}
N(2^{-2}\delta^2 c)^{1/4}  \right]$$

$$= (- \i z)^{-1/2} N(2^{-2}\delta^2
c)^{-1/4}f(\frac{-4}{c\delta^2z})$$

Observe that $(f\mid W(c))\mid W(c) = f.$
\item(The conjugation operator). The conjugation operator $H$ is
defined by $$(f\mid H)(z) = \overline{f(-\overline{z})}$$

\end{itemize}

\begin{lemma}\label{l:shift}

The operators $V(m), W(c)$ and $H$ take $M(\c, \psi)$ to $M(m\c,
\psi\epsilon_m) , M(\c, \overline{\psi}\epsilon_c)$ and $M(\c,
\overline{\psi})$ respectively. Further, if $f \in M(\c , \psi)$, we
have: \romanenum \begin{enumerate}

\item $(f \mid V(m))\mid T_{p^2} = (f \mid T_{p^2})\mid V(m) \quad$  when $p \nmid
m$

\item $(f\mid H)\mid T_{p^2} = (f \mid T_{p^2})\mid H$

\item $(f \mid W(c))\mid T_{p^2} = \psi_\c(p^2)(f \mid T_{p^2})\mid W(c)  \quad$ when $p \nmid
c$
\end{enumerate}
\end{lemma}

\begin{proof} The statements about $H$ are trivial while those about
$V(m)$ follow from (\cite{Shimura1}, Proposition 3.2) and
Proposition~\ref{p:Hecke} above.

We now prove the statements concerning $W(c)$. Let $W = \begin{pmatrix}
0 & -2\delta^{-1}  \\
2^{-1} \delta c & 0  \\
\end{pmatrix}$  and $\omega(z) = (- \i z)^{1/2}
N(2^{-2}\delta^2 c)^{1/4} $ ; by definition, $W(c)= [W, \omega(z)].$
Also, recall that $W_0
=\begin{pmatrix}
0 & -2\delta^{-1}  \\
2^{-1}\delta & 0  \\
\end{pmatrix}.$

Let us prove that $W(c)$ takes $M(\c, \psi)$ to $M(\c,
\overline{\psi}\epsilon_c)$. We need to show that
$$((f\mid W(c))
\parallel \gamma)(z) = \overline{\psi_\c(a_\gamma)}
(\epsilon_c)_\c(a_\gamma)(f \mid W(c))(z)$$ for any $\gamma \in
\Gamma_\c$.
 We write $$\Gamma = \left(
\begin{array}{ccc}
d_\gamma & -c_\gamma2^2\delta^{-2}c^{-1}  \\
-\delta^22^{-2}cb_\gamma & a_\gamma  \\
\end{array} \right) \in \Gamma_\c.$$  Using the fact that
$\Gamma^{-1}W \gamma = W$ and $(f \parallel \Gamma) = f$ we are
reduced to proving that
\begin{equation}\label{eqneedproveh}\frac{h(\Gamma , W z) )}{h(\gamma , z)} =
(\epsilon_c)_\c(a_\gamma) \frac{(- \i (\gamma z))^{1/2}}{(- \i
z)^{1/2}}.\end{equation}

Now put $\gamma' = V_0(c)\gamma V_0(c)^{-1}, z' =cz$. We note that $$h(\Gamma, Wz) = \theta (W_0\gamma'z')/ \theta(W_0 z').$$

Using the definition $h(G,z) = \theta(Gz)/\theta(z)$ for $G \in \D$ or $G = \begin{pmatrix}1/2&0\\0&2\end{pmatrix}W_0$ we have

 $$h(\Gamma, Wz) = \frac{h\left(\begin{pmatrix}1/2&0\\0&2\end{pmatrix}W_0 , V_0(1/4)\gamma'z'\right)h(V_0(c/4)\gamma V_0(c/4)^{-1}, cz/4)}{h\left(\begin{pmatrix}1/2&0\\
0&2\end{pmatrix}W_0, V_0(1/4)z'\right)}.$$ Use~\eqref{hwformula} on the factors $h\left(\begin{pmatrix}1/2&0\\0&2\end{pmatrix}W_0 , V_0(1/4)\gamma'z'\right)$, $h\left(\begin{pmatrix}1/2&0\\
0&2\end{pmatrix}W_0, V_0(1/4)z'\right)$ to get

 $$h(\Gamma , W z) = \frac{(- \i (\gamma z))^{1/2}}{(- \i
z)^{1/2}}h(V_0(c/4)\gamma V_0(c/4)^{-1}, cz/4) .$$ Now using~\eqref{e:hformula} we get $$ h(V_0(c/4)\gamma V_0(c/4)^{-1}, cz/4) = h(\gamma,z)(\epsilon_c)_\c(a_\gamma)$$ and this completes the proof of~\eqref{eqneedproveh}.

As for (iii), it follows from the following identities in $\G$,
which can be verified by explicit computation using the (partial)
automorphy property of $J_\Xi$.

\begin{enumerate}

\item Let $B =\begin{pmatrix}
1/p & 2b/(\delta p)  \\
0  & p \\
\end{pmatrix}$ where $b \in (R/p^2)^\times$. Suppose $b' \in (R/p^2)^\times$
be the element such that $bb'c \equiv -1$ mod $p^2$ and let $B'
=\begin{pmatrix}
1/p & 2b'/(\delta p)  \\
0  & p \\
\end{pmatrix}$.

Define $\gamma \in \Gamma_\c$ by $\gamma = \begin{pmatrix}
p^2 & -2b\delta^{-1}  \\
-2^{-1}\delta b' c  & (1+bb'c)/p^2 \\
\end{pmatrix} .$

Then we have
$$[\gamma, h(\gamma, z)]\ [B, J_\Xi(B, z)]\ [W, \omega(z)] = [W, \omega(z)]\ [B' ,
J_\Xi(B' , z)].$$

\item Let $C =\begin{pmatrix}
1 & 2h/(\delta p)  \\
0  & 1 \\
\end{pmatrix}$ where $h \in (R/p)^\times$. Suppose $h' \in (R/p)^\times$
be the element such that $hh'c \equiv -1$ mod $p$ and let $B'
=\begin{pmatrix}
1/p & 2h'\delta^{-1}  \\
0  & p \\
\end{pmatrix}$.\\

Define $\gamma ,\gamma' \in \Gamma_\c$ by $\gamma = \begin{pmatrix}
p & -2h\delta^{-1}  \\
-2^{-1}\delta h' c  & (1+hh'c)/p \\
\end{pmatrix}, \gamma' = \begin{pmatrix}
(1+hh'c)/p & -2h'\delta^{-1}  \\
-2^{-1}\delta h c  & p \\
\end{pmatrix} .$

Then we have
$$[\gamma, h(\gamma, z)]\ [C, J_\Xi(C, z)]\ [W, \omega(z)] = [W, \omega(z)]\ [B' , J_\Xi(B',z)].$$
$$[\gamma', h(\gamma', z)]\ [W, \omega(z)]\ [C, J_\Xi(C,z)] = (\epsilon_c)_\c(p)[B' , J_\Xi(B' , z)]\ [W,
\omega(z)].$$

\item Let $D =\begin{pmatrix}
p & 0  \\
0  & 1/p \\
\end{pmatrix}$ and $E
=\begin{pmatrix}
1/p & 0  \\
0  & p \\
\end{pmatrix}$.

Then we have
$$[D, N(p)^{-1/2}]\ [W, \omega(z)] = [W, \omega(z)]\ [E ,N(p)^{1/2}].$$
$$[W,\omega(z)]\ [D, N(p)^{-1/2}] =  [E , N(p)^{1/2}]\ [W,
\omega(z)].$$

\end{enumerate}

\end{proof}

Now, for a totally positive prime $p_0 \in R$  dividing $c/4$ let us write
$\Gamma_{\c /p_0}$ as a disjoint union of cosets modulo $\Gamma_\c$:
$$\Gamma_{\c /p_0} = \coprod_{\beta \in S}  \Gamma_\c \beta$$

We define the trace operator $S'(\psi) = S'(\psi, c, p_0)$ on $M(\c,
\psi)$ by $$(f\mid S'(\psi))(z) = \sum_{\beta \in S}
\psi(d_{\beta})(f
\parallel \beta) (z).$$ It is easy to see that this operator does not depend
on the choice of the $\beta$'s. Moreover if $r(\psi) \mid (\c/p_0)$,
$S'(\psi)$ takes $M(\c, \psi)$ to $M(\c/p_0 ,\psi)$ and if $f \in
M(\c / p_0 , \psi)$ then $f\mid S'(\psi) = uf$ where $u=\mid S\mid$.
A routine calculation similar to above also shows that $S'(\psi)$ commutes with
$T_{p^2}$ for $p  \nmid \c$.

We now define the operator $S(\psi) = S(\psi, c, p_0)$ on $M(\c,
\psi)$ by: $$S(\psi) = \frac{1}{u}N(p_0)^{1/4}
W(c)S'(\overline{\psi}\epsilon_c)W(c/p_0)
$$

\begin{lemma}

Let $p_0 \in R$ be a totally positive prime dividing $c/4$, such
that $r(\psi \epsilon_{p_0}) \mid (c/p_0)$. Then:
\begin{enumerate}

\item $S(\psi, c, p_0)$ maps $M(\c, \psi)$ into $M(\c/p_{0},
\psi\epsilon_{p_0})$

\item If $m$ is a totally positive element of $R$ that is prime to
$p_0$, and $f$ belongs to $M(\c, \psi)$, then
$$f\mid S(\psi, c, p_0)=f\mid S(\psi, mc, p_0).$$

\item $S(\psi)$ commutes with all the $T_{p^2}$

\item If $g \in M(\c/p_0, \psi\epsilon_{p_0})$, then $(g\mid V(p_0))\mid S(\psi, c,
p_0) = g$.

\item Let $p\in R$ be a totally positive prime such that $p \mid
(c/4)$, $p \neq p_0$ and $r(\psi \epsilon_{p}) \mid (c/p)$. If $g
\in M(\c /p, \psi \epsilon_{p})$, we have: $$(g\mid V(p))\mid
S(\psi, c, p_0) = (g\mid S(\psi \epsilon_{p}, c/p, p_0))\mid V(p).$$
\end{enumerate}
\end{lemma}

\begin{proof} The main ingredient for this proof was Lemma~\ref{l:shift}; otherwise the proof is identical to the proof of~\cite[Lemma 3]{Serre-Stark}.

(1) follows directly from Lemma~\ref{l:shift} and the comments above.

Now note that if $p_0 \nmid m$ and $ \gamma =\left(
\begin{array}{ccc}
a & b  \\
c & d  \\
\end{array} \right)$ runs over a set of representatives of
$\Gamma_{m\c} \backslash \Gamma_{m\c/p_0}$ then $ \gamma ' =\left(
\begin{array}{ccc}
a & bm  \\
c/m & d  \\
\end{array} \right)$ runs over a set of representatives of $\Gamma_{\c}\backslash \Gamma_{\c/p_0}$

To prove (2) we now only need to observe that $$W(mc)[\gamma,
h(\gamma,z)]W(mc/p_0) = [m I, 1]W(c)[\gamma', h(\gamma',
z)]W(c/p_0)$$

(3) follows from the commutativity of the Hecke operators with the
individual operators that make up $S(\psi)$.

As for (4) observe that $$(g\mid V(p_0))\mid W(c) =
N(p_0)^{-1/4}g\mid W(c/p_0).$$ The right side is invariant by
$\frac{1}{u}S'(\overline{\psi}\epsilon_c)$ and is sent to
$N(p_0)^{-1/4}g$ by $W(c/p_0)$.

Finally (5) follows from the following identities which can be checked by explicit computation:
$$\left[\begin{pmatrix}
p & 0  \\
0  & 1 \\
\end{pmatrix}, N(p)^{-1/4}\right]W(c) = [pI, 1]W(c/p),$$
$$W(c/p_0) = W(c/pp_0)\left[\begin{pmatrix}
p & 0  \\
0  & 1 \\
\end{pmatrix}, N(p)^{-1/4}\right].$$
.

\end{proof}

\begin{lemma}\label{l:g}
Let $m$ be a totally positive element of $R$ and $f\in M(\c, \psi)$.
Let $g(z) = f(z/m)$. Then $g \parallel \gamma =
\psi_\c(a_\gamma)\epsilon_m^\ast(a_\gamma)g$ for all $\gamma \in
\Gamma_{\c, m}$
\end{lemma}
\begin{proof}
Let $A=\begin{pmatrix}
1 & 0  \\
0  & p \\
\end{pmatrix}$ and let $\gamma' = A\gamma A^{-1}$. Note that $\gamma' \in \Gamma_\c$.

Then we have
$$[A, N(p)^{1/4}][\gamma , h(\gamma, z)] = [I,
\epsilon_m^\ast(a_\gamma)][\gamma', h(\gamma',z)][A, N(p)^{1/4}].$$
The result is now obtained by letting the above expression act on
$f$.
\end{proof}

Now, for any totally positive prime $p \mid (\c/4)$ we define the
operator $U(p) = U(p , \c)$ on $M(\c, \psi)$ by:
$$U(p) = N(p)^{-3/4} \sum_{j\in (R/p)}
\left[\begin{pmatrix}
1 & 2j\delta^{-1}  \\
0  & p \\
\end{pmatrix}, N(p)^{1/4}\right]
$$

\begin{lemma}\label{l:U}
$U(p)$ takes $M(\c, \psi)$ to $M(\c, \psi\epsilon_p)$ and if $$f(z)
= \sum_{\xi \in R} a(\xi)e(\xi z /2),$$ then $$(f\mid U(p))(z) =
\sum_{\xi \in R} a(p\xi)e(\xi z /2)$$
\end{lemma}

\begin{proof}
Let $A_j=\begin{pmatrix}
1 & 2j\delta^{-1}  \\
0  & 1 \\
\end{pmatrix}.$
 Put $g(z) = f(z/p).$ Note that $$f\mid \left[\begin{pmatrix}
1 & 2j\delta^{-1}  \\
0  & p \\
\end{pmatrix}, N(p)^{1/4}\right] = N(p)^{-1/4}g \parallel A_j$$

So for any $\gamma \in \Gamma_c$, $$f\mid \gamma = N(p)^{-1}
\sum_{j\in (R/p)}g\parallel (A_j\gamma).$$ But it is not hard to see
that $A_j$ varies over a set of right coset representatives of
$\Gamma_{\c,p}$ in $\Gamma_\c$; hence $A_j \gamma = \gamma'_i A_i$
with $\gamma'_i \in \Gamma_{\c,p}$ and distinct $j$ give rise to
distinct $i$. Note also that $a_{\gamma'_i} \equiv a_\gamma $ mod $
\c.$ Therefore \begin{align*} &N(p)^{-1}\sum_{j\in (R/p)}g\parallel
(A_j\gamma) \\& =N(p)^{-1} (\psi\epsilon_p)_\c(a_\gamma)\sum_{j\in
(R/p)}g\parallel A_i \\& =(\psi\epsilon_p)_\c(a_\gamma)f
\end{align*}

This proves that $U(p)$ takes $M(\c, \psi)$ to $M(\c,
\psi\epsilon_p)$.

As for the assertion about the Fourier
coefficients, note that \begin{align*}
(f\mid U(p))(z) &= N(p)^{-1}\sum_{j \in (R/p)}f(\frac{z + 2j \delta^{-1}}{p})\\ &=N(p)^{-1}\sum_{\xi \in R}a(\xi)e(\xi z/2p)\left(\sum_{j \in (R/p)} e(\frac{\xi j \delta^{-1}}{p})\right).
\end{align*}

The result now follows from the fact that $$\sum_{j \in (R/p)} e(\frac{\xi j \delta^{-1}}{p}) = \begin{cases} N(p) & \text{ if } p | \xi \\ 0 & \text{ otherwise }\end{cases}$$
\end{proof}

Finally, define the operator $K(p) = 1- U(p, p\c)V(p)$.

\begin{lemma}
If $f(z) = \sum_{\xi \in R} a(\xi)e(\xi z /2) \in M(\c , \psi)$ then
$f\mid K(p) \in M(\c p^2, \psi)$ and equals
$\sum_{(\xi,p)=1}a(\xi)e(\xi z /2)$. Further, if $p' \nmid p\c$ then
$T_{p'^2}$ and $K(p)$ commute.
\end{lemma}
\begin{proof}. This follows immediately from the above lemma and the
properties of $V(m)$ proved earlier. \end{proof}

\section{Newforms}\label{s:newforms}

\subsection{Definition of newforms and basic results}

Let $f \in M(\c, \psi)$ be an eigenvector of all but finitely many
$T_{p^2}$. We say that $f$ is an oldform if there exists a totally
positive prime $p$ dividing $\c /4$ such that one of the following
hold.

(a) $r(\psi)$ divides $(\c/p)$ and $f \in M(\c / p, \psi)$.

(b) $r(\psi \epsilon_p) \mid(\c/p)$ and $f = g\mid V(p)$ with $g\in
M(\c/p, \psi\epsilon_p)$.\\

We denote by $M^O(\c, \psi)$ the subspace of $M(\c , \psi)$ spanned
by oldforms. If $f \in M(\c, \psi)$ is an eigenvector of all but
finitely many $T_{p^2}$ and $f$ does not belong to $M^O(\c, \psi)$ ,
we say that $f$ is a newform of level $\c$.

The following two lemmas are proved exactly as in
\cite{Serre-Stark}. They are essentially formal consequences of all
the lemmas in the previous subsection.

\begin{lemma}
The symmetry operator and the conjugation operator take oldforms to
oldforms and newforms to newforms. \end{lemma}

\begin{lemma}
Let $h \in M^0(\c, \psi)$ be a non zero eigenform of all but
finitely many $T(p^2)$. Then there is a proper divisor $\c'$ of
$\c$, a character $\chi$ such that $r(\chi) \mid \c'$ and a newform
$g \in M(\c' , \psi)$ such that $h$ and $g$ have the same
eigenvalues for almost all $T(p^2)$.
\end{lemma}

We also have

\begin{lemma}
Let $p$ be a totally positive prime, and let $f(z) = \sum_{\xi \in
R} a(\xi)e(\xi z /2) \in M(\c , \psi)$ be non-zero and assume that
$a(\xi)=0$ for all $\xi$ not divisible by $p$. Then $p$ divides
$\c/4$, $r(\psi\epsilon_p)$ divides $\c/p$ and $f = g\mid V(p)$ with
$g\in M(\c/p, \psi\epsilon_p)$
\end{lemma}

\begin{proof}
Put $g(z) = f(z/p)$ and let $\c' = \c/p$ if $p \mid \c/4$ and $\c' =
\c$ otherwise. By Lemma~\ref{l:g} we have $$g\parallel \gamma =
(\psi\epsilon_p)_{p\c}(a_\gamma)g$$ for all $\gamma \in
\Gamma_{\c',p}$. Moreover, as $g$ has a Fourier expansion with
non-zero coefficients only in places corresponding to elements of
$R$, it follows that the above equation holds for $\gamma
=\begin{pmatrix}
1 & 2\delta^{-1}  \\
0  & 1 \\
\end{pmatrix}.$ By \cite[Lemma 3.4]{Shimura1}, the equation holds for
all $\gamma \in \Gamma_{\c'}$. Since $g$ is non-zero this implies
that $r(\psi\epsilon_p) | \c'$ which is possible only if $p$
divides $\c/4$. Thus $\c' = \c/p$ and hence $g\in M(\c/p,
\psi\epsilon_p).$

\end{proof}

The above lemmas allow us to derive our next theorem, which is the
main result that enables us to recognize oldforms. The proof of the
theorem is identical to that of Theorem 1 in \cite{Serre-Stark} and
will not be given here.
\begin{theorem}
Let $m$ be a totally positive element of $R$ and $f(z) = \sum_{\xi
\in R} a(\xi)e_{\mathbf{a}}(\xi z /2)$ be an element of $M(\c,
\psi)$ such that $a(\xi) = 0 $ for all $\xi$ with $(\xi, m)=1$.
Further assume that $f$ is an eigenform of all but finitely many
$T_{p'^2}$. Then $f \in M^0(\c , \psi)$. \end{theorem}

\subsection{Structure of newforms }
 Suppose $f(z) = \sum_{\xi \in R} a(\xi)e(\xi z /2) \in M(\c, \psi)$ is
 a newform. By theorem 3.1 there is a square free $t \in R$, unique
 up to multiplication by $U^2$, such that $a(\xi) = 0$ if $\xi/t$ is
 not a square.

The proofs of the next four lemmas are again identical to the
corresponding lemmas in \cite{Serre-Stark} and are omitted.
 \begin{lemma}\label{l:normalized}
 We have $t \in U^2$ and $a(1) \neq 0.$
 \end{lemma}

\begin{lemma}\label{l:scalar}
Let $g \in M(\c, \psi)$ be an eigenform of all but finitely many
$T(p^2)$, with the same eigenvalues as $f$. Then g is a scalar
multiple of $f$.
\end{lemma}

Because of Lemma~\ref{l:normalized}, we can divide by $a(1)$ and
henceforth assume that $f$ is normalized, i.e. $a(1) = 1$.

\begin{lemma}\label{l:eigenform}
Let $f \in M(\c, \psi)$ be a newform. Then $f$ is an eigenform for
every $T(p^2)$. Further, if $4p \mid \c$ , then the eigenvalue $c_p
= 0$.
\end{lemma}

\begin{lemma}\label{l:square}
The level $\c$ of the newform $f$ is a square and $f \mid W(c)$ is a
multiple of $f\mid H.$
\end{lemma}

\section{L-series and the proof of the main theorem}\label{main}
\subsection{The L-series}

Let $f(z) = \sum_{\xi \in R} a(\xi)e(\xi z /2)$ be an element of
$M(\c, \psi)$. For any ideal $I$, we define $a(I) = a(\xi)$ where
$\xi$ is any totally positive generator of $I$. Because $a(\xi u^2)
= \psi_\infty (u)a(\xi)$ for any unit $u$ by (\cite{Shimura1}, Prop.
5.4), it follows that if we assume that $\psi_\infty$ is trivial,
then $a(I)$ is well-defined. In that case we define the $L$-series
$L(s,f)$ by
$$L(s,f) = \sum \frac{a(I)}{N(I)^s}$$ where the sum is taken over
all non-zero ideals of $R$.

\begin{theorem}\label{t:lseries}

Suppose $f\in M(\c, \psi)$ where $\psi_\infty$ is trivial and assume
that $\c$ is the square of an ideal. Then $L(s,f)$ can be
analytically continued to an entire function (with the exception of
a simple pole at $s = 1/2$ if $f$ is not a cusp form). Moreover, if
$$\Lambda(s ,f) = (2\pi)^{-ns} \Gamma(s)^n N(\delta)^s
N(\c)^{(s/2)}L(s,f)$$ then the following relation holds
$$\Lambda(s,f) =  \Lambda(1/2 -s , g)$$ where $g = f\mid W(c)$ .

\end{theorem}
\begin{proof}
Let $$f(z) = \sum_{\xi \in R} a(\xi)e(\xi z /2)$$ and
$$g(z)= \sum_{\xi \in R} b(\xi)e(\xi z /2).$$ Also let $\c =
(c_1)^2$ where $c_1$ is totally positive (we can do this because $F$
has narrow class number one).

Put $f_1 = f - a(0)$ , $g_1 = g - b(0)$. Recall that
$F_{\infty}^\circ$ denotes $\prod_{v \in \infty } F_v^{+}$ which can
be naturally identified with $(\R^+)^n$. Thus there is an action of
$U^2$ on $(\R^+)^n$ and for later purposes it is important to note
that this action preserves the norm of an element. Now consider the
coset space $(\R^+)^n/U^2$. Define the integral
\begin{equation}\label{e:integral} \Phi(s) =\int_{(\R^+)^n/U^2}f_1(\frac{2\i
y}{c_1\delta})\prod_{j=1}^{n}y_j^s\frac{dy_j}{y_j}.\end{equation}

We first observe that this integral is convergent for all $s$ with
$Re(s)>1/2$. Indeed, by the unit theorem, we may choose the
fundamental domain $F_{\infty}^\circ/U^2$ such that the ratios
$y_i/y_j$ are all bounded, and hence all the $y_j$ go to zero or
infinity together. As the $y_j \rightarrow \infty$ the rapid decay
of $f_1$ assures convergence. As they go to zero, we use the
following equation, which follows easily from $g = f\mid W(c)$:

\begin{equation}\label{e:fg}
 f_1(\frac{2\i }{c_1\delta y}) = \prod y_j^{1/2}g_1(\frac{2\i
y}{c_1\delta}) \quad + b(0)\prod y_j^{1/2} \quad - a(0)
\end{equation}
to obtain the same result.

Now, write the right side of (\ref{e:integral}) as
\begin{align*}
& \quad \quad \sum_{I\neq 0} a(I) \int_{(\R^+)^n/U^2} \sum_{(\alpha)
= I, \alpha
>>0}e^{-2\pi tr(\alpha y/(c_1\delta))}
\prod_{j=1}^{n}y_j^s\frac{dy_j}{y_j}.\\
& =\sum_{\alpha \in R^+/U^2}\sum_{\epsilon \in
U^2}a(\alpha\epsilon)\int_{(\R^+)^n/U^2}e^{-2\pi tr(\alpha\epsilon
y/(c_1\delta))}\prod_{j=1}^{n}y_j^s\frac{dy_j}{y_j}.\\
& =\sum_{\alpha \in R^+/U^2}a(\alpha)\int_{(\R^+)^n/U^2}e^{-2\pi
 tr(\alpha y/(c_1\delta))}\prod_{j=1}^{n}y_j^s\frac{dy_j}{y_j}.\\
 &
=\sum_{\alpha \in
R^+/U^2}a(\alpha)\prod_{j=1}^n[(2\pi)^{-s}(c_1^{(j)}\delta^{(j)}/\alpha^{(j)})^s\Gamma(s)]\\
&=(2\pi)^{-ns}\Gamma(s)^n N(\delta)^s N(\c)^{s/2}L(s,f)\\
&=\Lambda(s,f).
\end{align*}

On the other hand, we have,
\begin{align*}
& \quad \int_{y\in (\R^+)^n/U^2, N(y)<1}f_1(\frac{2\i
y}{c_1\delta})\prod_{j=1}^{n}y_j^s\frac{dy_j}{y_j}\\
& =\int_{y\in (\R^+)^n/U^2, N(y)<1}(f(\frac{2\i
y}{c_1\delta})-a_0)\prod_{j=1}^{n}y_j^s\frac{dy_j}{y_j}\\
&=-a_0\int_{y\in (\R^+)^n/U^2,
N(y)<1}\prod_{j=1}^{n}y_j^s\frac{dy_j}{y_j} + \int_{y\in
(\R^+)^n/U^2,
N(y)>1}f(\frac{2\i}{yc_1\delta})\prod_{j=1}^{n}y_j^{-s}\frac{dy_j}{y_j}\\
&= -a_0\int_{y\in (\R^+)^n/U^2,
N(y)<1}\prod_{j=1}^{n}y_j^s\frac{dy_j}{y_j} + \int_{y\in
(\R^+)^n/U^2, N(y)>1}(g(\frac{2\i y}{c_1\delta}) -
b_0)\prod_{j=1}^{n}y_j^{1/2 - s}\frac{dy_j}{y_j} \\
&\quad +b_0\int_{y\in (\R^+)^n/U^2,
N(y)<1}\prod_{j=1}^{n}y_j^{s-1/2}\frac{dy_j}{y_j}\\
&= -\frac{a_0C}{s} - \quad \frac{b_0C}{(1/2 -s)} +\quad \int_{y\in
(\R^+)^n/U^2, N(y)>1}g_1(\frac{2\i
y}{c_1\delta})\prod_{j=1}^{n}y_j^{(1/2 - s)}\frac{dy_j}{y_j}.
\end{align*}

for some constant $C$. Note that in the last step, we have used the
Dirichlet unit theorem.

Hence we have shown that
\begin{align*}
 (5) \quad \quad \quad \quad \quad \quad \Phi(s) + \frac{a_0C}{s}
+ \frac{b_0C}{(1/2 -s)} = &\int_{y\in (\R^+)^n/U^2,
N(y)>1}[f_1(\frac{2\i
y}{c_1\delta})\prod_{j=1}^{n}y_j^s\frac{dy_j}{y_j}\\
& \quad +  g_1(\frac{2\i y}{c_1\delta})\prod_{j=1}^{n}y_j^{(1/2 -
s)}\frac{dy_j}{y_j} ]
\end{align*}

The right side consists of integrals over regions that are bounded
away from 0 (by the Dirichlet unit theorem) and hence the rapid
decay of $f_1$ and $g_1$ near infinity imply that these integrals
converge for all s. This proves that $\Phi(s)$ is a meromorphic
function with simple poles at $0, 1/2$ if $f$ is not a cusp form. As
a corollary, we obtain that $L(s,f)$ can be analytically continued
to the entire complex plane (with a simple pole at $1/2$ if $f$ is
not a cusp form).

To see the functional equation just exchange the roles of $f$ and
$g$ in (5).

\end{proof}

We call a prime ideal $\p$ of $R$ \emph{non-split} if $\p$ is the
unique prime ideal of $R$ that lies above $\p \cap \Z$. We call an
ideal $I$ of $R$ non-split if all its prime divisors are non-split.

\begin{theorem}\label{t:newtheta}
Suppose $\c$ is non-split. Let $f$ be a normalized newform in $M(\c,
\psi)$ with $\psi_\infty$ trivial. Then $\c = 4r(\psi)^2$ and $f =
\frac{1}{2}\theta_{\psi}$
\end{theorem}

\begin{proof}. By Theorem~\ref{t:serre}, Lemma~\ref{l:normalized} and Lemma~\ref{l:eigenform}, we have the
product decomposition
$$L(s,f) = \prod_{\p \mid \c}\left(1-\frac{c_p}{N(\p)^{2s}}\right)^{-1}
\prod_{\p \nmid
\c}\left(1-\frac{\psi^\ast(\p)}{N(\p)^{2s}}\right)^{-1}$$
Furthermore by Lemma~\ref{l:square} and Theorem~\ref{t:lseries} we
have
$$(2\pi)^{-ns} \Gamma(s)^n L(s,f) = C_1(2\pi)^{-n(1/2 -s)}
(\Gamma(1/2 - s))^n N(c\delta^2)^{1/2 -s}L(1/2 -s,Hf)$$ for some
constant $C_1$.

Consider, on the other hand, the function $L(2s, \psi)$ defined by
$$L(2s, \psi) = \frac{\psi^\ast(I)}{N(I)^{2s}} = \prod_{\p \nmid r(\psi)}\left(1 -
\frac{\psi^\ast(\p)}{N(\p)^{2s}}\right)^{-1}$$

Then, from (\cite{Bump}, p. 78-79) we know that $$(2\pi)^{-ns}
\Gamma(s)^n L(2s, \psi) = C_2 (2\pi)^{-n(1/2 -s)}N(4 r(\psi)^2
\delta^2)^{1/2 -s}(\Gamma(1/2 -s))^n L(1- 2s ,\overline{ \psi})$$

Dividing these equations we have $$\prod_{\p \in S}(\frac{1-c_p
N(\p)^{-2s}}{1-\psi^\ast(\p)N(\p)^{-2s}}) = C_3 N(\c /
4r(\psi)^2)^{-(1/2-s)}\prod_{\p \in S}(\frac{1-\overline{c_p}
N(\p)^{2s - 1}}{1-\overline{\psi^\ast(\p)}N(\p)^{2s - 1}})$$ where
$S$ is the set of prime ideals $\p$ for which $c_p \neq
\psi^\ast(\p)$,  $\p \mid \c$.

If, for some $\p \in S$, we have $\psi^\ast(\p) \neq 0$, then the
left side of the above equation has an infinity of poles on the line
$Re(s) = 0$, only finitely many of which can appear on the right
side. This can be seen as follows: if $\p$ is a prime in S then by
assumption, it is the only prime with that norm and so we can find
infinitely many $s$ such that $\psi^\ast(\p)N(\p)^{-2s} = 1$ but
none of the expressions $c_{p'} N(\p')^{-2s}, \overline{c_{p'}}
N(\p')^{2s - 1}, \overline{\psi^\ast(\p')}N(\p')^{2s - 1}$ equals 1
for any $\p' \in S, \p' \neq \p.$

Hence $\p \in S$ implies $\psi^\ast(\p) = 0$, (in other words $\p
\mid r(\psi)$) and hence $c_p \neq 0$ since $c_p \neq
\psi^\ast(\p)$. But $c_p = 0$ if $4\p \mid \c$. It follows that
either $S$ is empty or consists of the unique prime  that lies above
$2$. Meanwhile, the equation simplifies to
$$\prod_{\p \in S}(1-c_p N(\p)^{-2s}) = C_4 N(\c \m^2 /
4r(\psi)^2)^s\prod_{\p \in S}(1-c_p' N(\p)^{-2s})$$ where $c_p' =
N(\p)/\overline{c_p}$ and $\m = \prod_{\p \in S}\p$.

We claim that $S$ is empty. Suppose not, then $S = \{ \p \}$ where
$\p$ is the unique prime above $2$. Then, if $c_p \neq c_{p'}$ we
can find a zero of the left side of the above identity that is not a
zero of the right side. Hence we must have $c_p = c_{p'}$. This
implies that $\mid c_p^2  \mid= N(\p)$. But that contradicts
Corollary~\ref{cbound2}

Thus $S$ is empty and we have $\c = 4r(\psi)^2$, $L(2s , \psi) = L(s
, f)$. This implies that for any non zero ideal $I$ which has a
common factor with $\c$ we have that $a(I) = \psi^\ast(L)$ if $I =
L^2$ for some ideal $L$ coprime to $r(\psi)$ and $a(I) = 0$ in all
other cases. This, coupled with Theorem~\ref{t:serre} shows that $f$
and $\frac{1}{2}\theta_\psi$ have the Fourier coefficients at $\xi$
for all $\xi \neq 0$; hence they also have the same constant
coefficient. Thus $f = \frac{1}{2}\theta_\psi.$  \end{proof}

\subsection{Proof of Theorem~\ref{t:main}}

\begin{proof}. We break the proof into two parts:

\begin{case}The $\theta_{\psi, t}$ are linearly independent. \end{case}

Since $t$ and $\psi$ determine $\chi$, each $t$ occurs as the second
entry of at most one $(\chi, t)$ in $\Omega(\c, \psi)$. Suppose we
have
$$\lambda_1\theta_{\psi_1, t_1} + \lambda_2\theta_{\psi_2, t_2}
+...+\lambda_m\theta_{\psi_m, t_m} = 0$$ with the number of primes
in the prime decomposition of $t_1$ being less than or equal to that
for the other $t_i$ and $\lambda_i \neq 0$ for each $i$. Then the
coefficient at place $t_1$ is $2\lambda_1$ for $\theta_{\psi_1,
t_1}$ and 0 for the others, thus showing that $\lambda_1 = 0$, a
contradiction.

\begin{case} The $\theta_{\psi, t}$
span $M(\c, \psi)$. \end{case}

We use induction on the number of (not necessarily distinct)prime
factors of $\c$. By Lemma~\ref{l:basis}, it suffices to show that any
eigenform $f$ of all the $T_{p^2}, p\nmid \c$ is a linear
combination of the $\theta_{\chi,t}$ with $(\chi, t) \in \Omega(\c,
\psi)$. If $f$ is a newform, this follows from
Theorem~\ref{t:newtheta}. If not, we may assume $f$ is an oldform.
Now we have two cases.

In the first case, $r(\psi)$ divides $\c/ p$ and $f\in M(\c/p,
\psi)$. Since $\c/ p$ is also non-split the induction hypothesis
shows that $f$ is a linear combination of the $\theta_{\chi, t}$
with $(\chi, t)$ in $\Omega(\c/p, \psi)$ and hence in $\Omega(\c ,
\psi)$.

In the second case, $r(\psi\epsilon_p)$ divides $\c/p$ and $f =
V(p)g$ with $g \in M(\c/p, \psi\epsilon_p)$. Because $\c /p$ is
non-split, and $(\psi \epsilon_p)$ is totally even because $\psi$ is,
the induction hypothesis shows that $g$ is a linear combination of
the $\theta_{\chi,t}$ with $(\chi,t) \in \Omega(\c/p,
\psi\epsilon_p)$ and hence $f$ is a linear combination of the
$\theta_{\chi, tp}$, with $(\chi, tp) \in \Omega(\c , \psi)$. This
completes the proof.

\end{proof}

\subsection{Examples}

In this section we specialize to the case $F= \Q(\sqrt{2})$, and $\c
= (\sqrt{2})^n$. For brevity, let $q = \sqrt{2}$. Note that the ring
of integers $R$ is simply $\Z(\sqrt{2})$ and the unit group is $<-1>
\times <1+q>$. Note also that $F$ has narrow class number 1, and the
prime $2$ is ramified in $F$, which allows us to apply the theorems
of the last section with $\c = (q)^n$.

The theorem will apply to any Hecke character $\psi$ of $F$ with
$\psi_\infty$ trivial and such that $r(\psi)$ divides $q^n$. For
simplicity we only find the quadratic (of order 2) Hecke characters
if this type, and give the explicit bases for each of the
corresponding spaces of modular forms. It suffices to find the
quadratic Dirichlet characters mod $q^n$ that are trivial on units.
For that we need to analyze the structure of the groups
$(R/q^n)^\times$.

\begin{proposition}\label{p:anal}
Let $U_n$ denote the multiplicative group $(R/q^n)^\times$. Then, if
$n \leq 4$, $U_n$ is generated by the units of $R$ and hence there
is no nontrivial even Hecke character with conductor dividing $q^n$.
On the other hand, if $n>4$, the following hold: \begin{enumerate}

\item $U_n$ is the direct sum of the cyclic groups generated by $(1+q)
, (-1)$ and $(3+4q)$.

\item $(1+q)$ has order $2^{\lfloor \frac{n}{2}\rfloor}$ while $3+4q$
has order $2^{\lfloor \frac{n-3}{2}\rfloor}$ in the group
$(R/q^n)^\times$.

\item Let $k=\lfloor \frac{n}{2}\rfloor$, $l= \lfloor
\frac{n-3}{2}\rfloor$. Then $U_n$ is isomorphic to $(\Z/2^k\Z)
\oplus (\Z/2^l\Z) \oplus (\Z/2\Z)$

\end{enumerate}
\end{proposition}

\begin{proof}. The case $n \leq 4$ can be checked easily by hand.

For the case $n \geq 5$, first observe that the cardinality of $U_n$
is $2^{n-1}$. This follows from the fact that the elements of
$(R/q^n)$ can be written as $(a+bq)$ with $a \in (\Z/2^k\Z)$ and $b
\in (\Z/2^{l+1}\Z)$ with $k,l$ as in part (3) of the Proposition,
and such an element is invertible iff $a$ is odd.

We first prove (2). The same method is used to calculate the orders
of $3+4q$ and $1+q$; the idea is to write $(a+bq)^{2^{k+1}} - 1 =
(a+bq)^{2^k} - 1)(a+bq)^{2^k} + 1).$ If for some $n>2$, $q^n$
\emph{exactly} divides $(a+bq)^{2k} - 1)$ then $q^{n+2}$ exactly
divides $(a+bq)^{2^{k+1}} - 1 $. Since $q^5$ exactly divides
$(1+q)^4 -1$ and $q^6$ exactly divides $(3+4q)^2 -1$, the result
follows.

Now, it is easy to see that the subgroups generated by $(1+q)$ and
$(3+4q)$ have trivial intersection, and further, that $(-1)$ does
not lie in the subgroup generated by these two elements. Thus (a)
follows by comparing cardinalities, and clearly (3) is a direct
consequence of (1) and (2).

\end{proof}

\begin{corollary}
Let $\phi$ denote the Hecke character $\epsilon_u$ where $u=2+q$.
Then $r(\phi) = (q^5)$ and $\phi$ is the unique non-trivial Hecke
character with $\psi_\infty$ trivial that is quadratic (of order 2)
and whose conductor divides $q^n$.
\end{corollary}
\begin{proof}
Observe that any Dirichlet character mod($q^n$) that is
trivial on the units of $R$ must be, by the previous proposition, a
character on the group generated by $3+4q$. Furthermore, if this
character is quadratic, it must be either the trivial character or
the character that takes value $-1$ on $(3+4q)$; it is not hard to
see (since $(3+4q)$ is inert in the extension $F(\sqrt{u})/F$) that
this corresponds to the Hecke character $\epsilon_u$. \end{proof}
This leads us to the following theorem.

\begin{theorem}
Let $n \geq 5$, $\c = (q^n)$ , $u = 2+q$. Let $\phi$ denote the
Hecke character $\epsilon_u$ and $\mathbf{1}$ denote the trivial
character. Then: \begin{enumerate}

\item A basis for $M(\c, \mathbf{1})$ comprises of the functions $\{
\theta_{\mathbf{1}, 2^k} , 0 \leq k \leq \lfloor \frac{n-4}{2}
\rfloor$; $\theta_{\phi, 2^ku} , 0 \leq k \leq \lfloor
\frac{n-15}{2} \rfloor \}$. Thus the dimension of the space $M(\c,
\mathbf{1})$ is $( \lfloor \frac{n-2}{2} \rfloor + max\{\lfloor
\frac{n-13}{2} \rfloor, 0\} )$.

\item A basis for $M(\c, \phi)$ comprises of the functions $\{
\theta_{\mathbf{1}, 2^ku} , 0 \leq k \leq \lfloor \frac{n-5}{2}
\rfloor$; $\theta_{\phi, 2^k} , 0 \leq k \leq \lfloor \frac{n-14}{2}
\rfloor \}$. Thus the dimension of the space $M(\c, \mathbf{1})$ is
$( \lfloor \frac{n-3}{2} \rfloor + max\{\lfloor \frac{n-12}{2}
\rfloor, 0\} )$.

\end{enumerate}
\end{theorem}

\begin{proof}. This follows from Theorem~\ref{t:main} and the above Corollary.
\end{proof}

\section{Potential applications }\label{s:app}
In this section we put our work in context by mentioning a couple of potential applications which we hope to take up elsewhere.
\subsection{The Congruence number problem}
An ancient Diophantine problem (the so-called congruence number problem) asks for a good criterion to determine whether an integer is the area of a right angled triangle with rational sides. Such integers are referred to as congruent numbers.
This was solved by Tunnell~\cite{Tunnell}. Tunnell's work begins with the observation that $n$ is congruent if and only if the rank of the elliptic curve $E = y^2 = x^3 - n^2x$ over $\Q$ is non-zero. This is easy to prove by elementary number theory. Now, by the Birch--Swinnerton-Dyer conjecture (one direction of which is known in this case by the work of Coates--Wiles) the above condition is equivalent to the value of the $L$-function $L(E,s)$ at 1 (the central value) being equal to 0. However, it is not hard to show that $L(E,s)$ equals $L(\phi \otimes \epsilon_n,s)$ where $\phi$ is the unique normalized newform of weight 2, level 32 and trivial character while $\epsilon_n$ is the quadratic character associated to $\Q(\sqrt{n})$.    By work of Waldspurger the value $L(\phi\otimes\epsilon_n,1)$ is related to the value $c_n ^2$ where $c_n$ is the $n$'th Fourier coefficient of the weight $3/2$ modular form that maps to $\phi$ under the Shimura correspondence.

Tunnell's main contribution to the problem was to find explicitly the weight $3/2$ form above. Using the Serre-Stark theorem, he was able to write this form as a product of an explicit theta-series and a standard weight 1 form. As a result, it was possible to express $c_n^2$ and consequently the vanishing condition on $L(E,1)$ in a simple combinatorial form.

One may ask the same question over our totally real number field $F$. We call an element $\alpha \in R$, $F$-congruent if there exist  positive $X, Y, Z \in F$ such that $X^2 + Y^2 = Z^2$ and $XY = 2 \alpha$, with possibly a signature restriction. In the case of real quadratic fields, things work out nicely, though with a slight modification~\cite{Achimescu}. Thus we hope that one can resolve the congruent number problem over $F$ in a manner similar to what was achieved by Tunnell over $\Q$. One of the crucial points is the construction of an appropriate weight $3/2$ Hilbert modular form; we hope to achieve this by using our basis of weight $1/2$ forms and multiplying it by a appropriate weight $1$ Hilbert modular form.

\subsection{Construction of interesting weight $1$ forms} There are not many explicit examples that illustrate the conjectural correspondence between Galois representations and weight $1$ forms. Buhler~\cite{Buhler} was able to construct a (classical) modular form of level 800 that corresponds to an icosahedral Galois representation. We believe that our main theorem may be useful in constructing interesting (that is, not a base change and non-dihedral) Hilbert modular forms of a specified level whose $L$- function matches a (possibly icosahedral) Galois representation of that level.

Given a Hilbert modular form $f$ of weight $1$ and two forms $g_1, g_2$ of weight $1/2$, $fg_1g_2$ is a form of weight 2. This suggests the following procedure. We fix a form $F$ of weight $2$ and level $N$ and consider the functions $F/(g_1g_2)$ where $(g_1,g_2)$ varies over pairs of basis forms of weight $1/2$, as given by our main theorem. It seems likely that this method will lead to the construction of explicit, interesting examples of weight 1 Hilbert modular forms that correspond to Galois representations.

\bibliography{hilbert}

\def\cprime{$'$}
\begin{thebibliography}{10}

\bibitem{Achimescu}
S~Achimescu.
\newblock {\em Hilbert Modular forms of weight $1/2$}.
\newblock PhD thesis, Caltech, 2004.

\bibitem{Buhler}
Joe Buhler.
\newblock An icosahedral modular form of weight one.
\newblock In {\em Modular functions of one variable, V (Proc. Second Internat.
  Conf., Univ. Bonn, Bonn, 1976)}, pages 289--294. Lecture Notes in Math., Vol.
  601. Springer, Berlin, 1977.

\bibitem{Bump}
Daniel Bump.
\newblock {\em Automorphic forms and representations}, volume~55 of {\em
  Cambridge Studies in Advanced Mathematics}.
\newblock Cambridge University Press, Cambridge, 1997.

\bibitem{Garrett}
Paul~B. Garrett.
\newblock {\em Holomorphic {H}ilbert modular forms}.
\newblock The Wadsworth \& Brooks/Cole Mathematics Series. Wadsworth \&
  Brooks/Cole Advanced Books \& Software, Pacific Grove, CA, 1990.

\bibitem{gelpia}
Stephen Gelbart and I.~I. Piatetski-Shapiro.
\newblock Distinguished representations and modular forms of half-integral
  weight.
\newblock {\em Invent. Math.}, 59(2):145--188, 1980.

\bibitem{Serre-Stark}
J.-P. Serre and H.~M. Stark.
\newblock Modular forms of weight {$1/2$}.
\newblock In {\em Modular functions of one variable, VI (Proc. Second Internat.
  Conf., Univ. Bonn, Bonn, 1976)}, pages 27--67. Lecture Notes in Math., Vol.
  627. Springer, Berlin, 1977.

\bibitem{Shimura3}
Goro Shimura.
\newblock On modular forms of half integral weight.
\newblock {\em Ann. of Math. (2)}, 97:440--481, 1973.

\bibitem{Shimura2}
Goro Shimura.
\newblock On {E}isenstein series of half-integral weight.
\newblock {\em Duke Math. J.}, 52(2):281--314, 1985.

\bibitem{Shimura1}
Goro Shimura.
\newblock On {H}ilbert modular forms of half-integral weight.
\newblock {\em Duke Math. J.}, 55(4):765--838, 1987.

\bibitem{Shimura4}
Goro Shimura.
\newblock {\em Euler products and {E}isenstein series}, volume~93 of {\em CBMS
  Regional Conference Series in Mathematics}.
\newblock Published for the Conference Board of the Mathematical Sciences,
  Washington, DC, 1997.

\bibitem{Tunnell}
J.~B. Tunnell.
\newblock A classical {D}iophantine problem and modular forms of weight
  {$3/2$}.
\newblock {\em Invent. Math.}, 72(2):323--334, 1983.

\end{thebibliography}

\end{document}